\documentclass[12pt]{article}
\usepackage{amssymb}
\usepackage{amsmath}
\usepackage{amssymb}
\usepackage{amsthm}
\usepackage{graphicx}
\usepackage{amsthm}
\usepackage{graphicx}
\usepackage[dvipsnames]{xcolor}
\usepackage{graphics}
\usepackage{amsthm}

\input pictex.tex

\theoremstyle{definition}

\input epsf

\begin{document}

\date{}
\author{Andr\'es Navas}

\title{Parametric proofs of the Pythagorean theorem \\
via ziggurats and pyramids}
\maketitle

\noindent{\bf Abstract.}  We present two parametric families of area-rearrangement proofs of the Pytha- gorean 
theorem. Our contribution is primarily constructive and expository: we organize the geometry through a single
 angular parameter. For certain appropriate parameters, we recover several known motifs as special instances. 
 For less regular parameters, the argument becomes more complicated, requiring the use of trigonometric 
 identities, a topic we also explore in further detail.

\vspace{0.2cm}


\vspace{0.2cm}

\noindent{\bf MCS 2020:} 
01-02, 
51-03, 
97G30, 
97G99, 
97U99. 


\vspace{1cm}

\noindent{\bf \Large{Introduction}}

\vspace{0.5cm}

In recent years, there has been some interest in (re)discovering and collecting different proofs of the Pythagorean theorem, 
much beyond the classical compilation \cite{Lo}. Here we propose two new proofs, which come from interesting geometric 
configurations that, apparently, have remained unnoticed in their full generality. In fact, these are two parametric families of proofs, 
each of which depends on an angle involved in the corresponding construction (cf. \cite{Ho}). 

For the first family, the parameter is the angle of inclination of a trapezium that, as 
a little joke, we call a {\em ziggurat}. In concrete terms, by a ziggurat of side $\ell$ and 
angle $60^{\mathrm{o}} \leq \theta < 180^{\mathrm{o}}$ (an $(\ell,\theta)$-ziggurat, for short) we mean the isosceles 
trapezium having a base and the two ``nonparallel'' sides of length $\ell$, with both angles on the basis equal to $\theta$. 

\vspace{-0.2cm}

\begin{center}
\includegraphics[width=6.5cm]{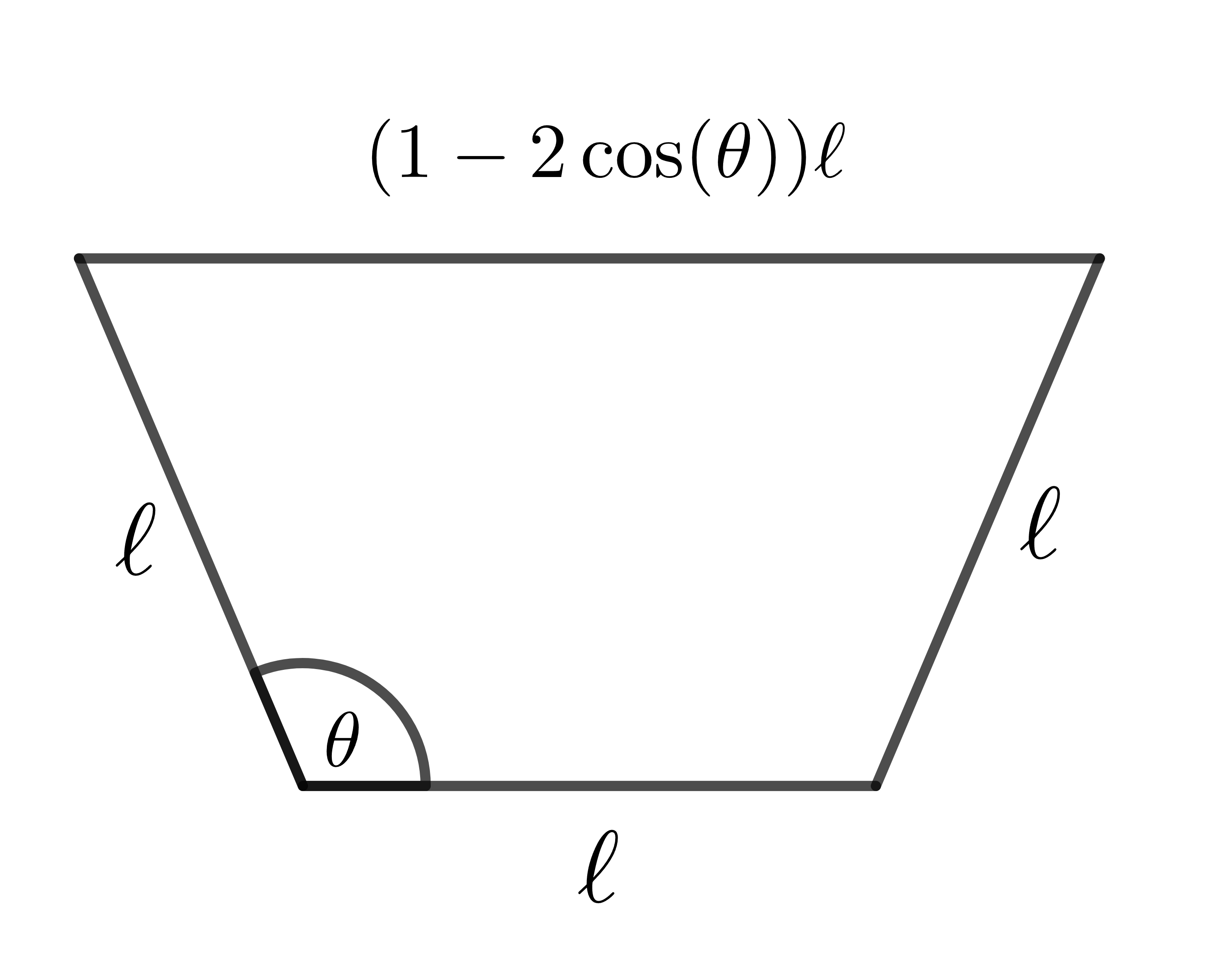}
\end{center}

\vspace{0.2cm}

Note that, for $\theta = 90^{\mathrm{o}}$, this becomes a square. Also, for $\theta = 60^{\mathrm{o}}$, this degenerates to an equilateral 
triangle. The first aim of this note is to give a direct proof (with no use of the Pythagorean theorem) of the following result. Let us emphasize 
that the new contribution here is not the theorem below itself but the proof argument; indeed, as pointed out in \cite{PS}, the result is known since 
the time of Euclid (in an even more general form)\footnote{More concretely, Proposition 31 of Book VI of Euclid's {\em Elements} 
states: ``In right-angled triangles, the figure on the side opposite to the right angle equals the sum of the 
similar and similarly described figures on the sides containing the right angle'', and this for any figure~! (of course, we know today 
that nonmeasurable figures are forbidden).}  as a direct consequence of the Pythagorean theorem. 

\vspace{0.4cm}

\noindent{\bf Theorem A.}
Let $\triangle ABC$ be a triangle with a right angle at $C$ and sides $a,b$ and $c$ (opposed to $A,B$ and $C$, respectively). 
For each $60^{o} \leq \theta \leq 135^{o}$, the sum of the areas of the $(a,\theta)$ and $(b,\theta)$-ziggurats equals the area 
of the $(c,\theta)$-ziggurat.

\begin{center}
\includegraphics[width=10cm]{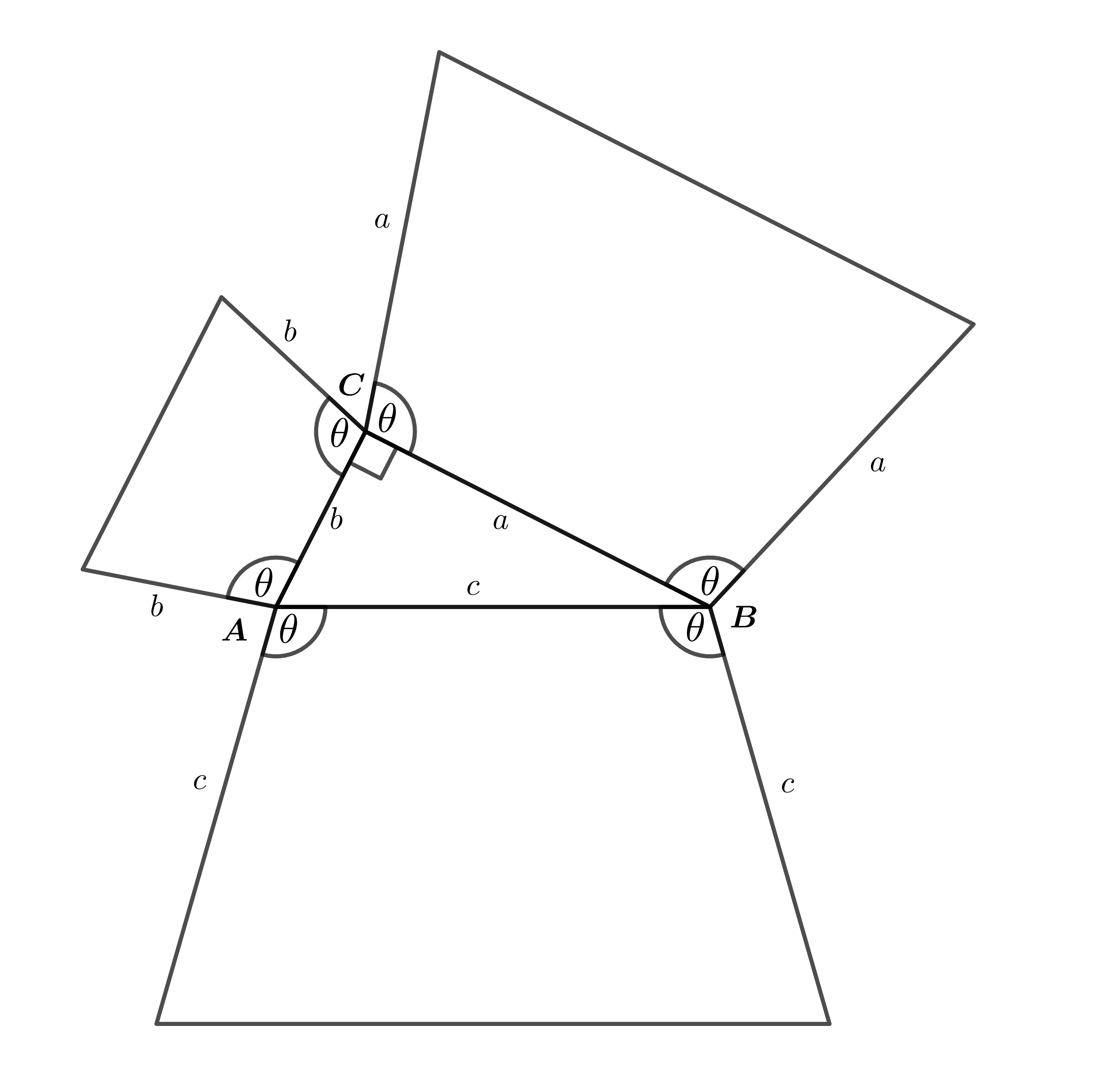}
\end{center}

For $\theta = 90^{\mathrm{o}}$, this corresponds exactly to the Pythagorean theorem. 
However, for each angle $\theta$, the corresponding claim is also equivalent to it.\footnote{This idea appears in many sources in 
the literature; in particular, it is widely discussed in \cite{PS} (unfortunately, we were unaware of this reference while writing \cite{Na1}, 
where the idea also plays a pivotal role).} Indeed, the area of the $(\ell,\theta)$-ziggurat equals some constant $C(\theta)$ times the 
area of the square of side $\ell$, so that dividing all the terms in the equality involving the areas of the ziggurats by $C(\theta)$, 
we recover the equality involving the squares. 
Here, the value of the (positive) constant $C(\theta)$ is irrelevant, yet it can explicitly (and easily) be computed:
 $C(\theta) = \sin (\theta) \, (1 - \cos (\theta)).$

If we consider ziggurats of an angle $90^{\mathrm{o}} < \theta \leq 135^{\mathrm{o}}$ and we extend its non parallel sides until they meet, 
we obtain isosceles triangles over the sides of our original triangle $\triangle ABC$, each having two angles equal to $\theta$ at the basis. 
The area of such a triangle equals some universal constant $D(\theta)$ times that of the corresponding ziggurat. (Again, the precise 
value of $D(\theta)$ is irrelevant, yet it is easily seen to be equal to $1 / (4 \cos (\theta) (1-\cos (\theta)))$.) Thus, Theorem A implies 
the result below, which, of course, is also a consequence of the Pythagorean theorem. To simplify the statement, as another little joke, 
we call $(\ell,\theta)${\em -pyramid} an isosceles triangle with basis $\ell$ and angles $\theta$ at the basis.\footnote{The reader may 
complain in that ziggurats and pyramids are terms that are more appropriate for 3D figures. However, the joke in the terminology 
is just a stylistic way to recall that the Pythagorean theorem was somehow known to civilizations prior to the Greeks, notably in ancient 
Mesopotamia and Egypt.} 

\vspace{0.4cm}

\noindent{\bf Theorem B.}
Let $\triangle ABC$ be a triangle with a right angle at $C$ and sides $a,b$ and $c$ (opposed to $A,B$ and $C$, respectively). For each 
$45^{\mathrm{o}} \leq \theta < 90^{\mathrm{o}}$, the sum of the areas of the $(a,\theta)$ and $(b,\theta)$-pyramids equals the area 
of the $(c,\theta)$-pyramid. 

\vspace{0.55cm}

Again, our contribution here is the proof argument, which is given in the last section of this note. It is worth mentioning that 
this argument is direct, that is, it does not pass through ziggurats and, of course, avoids the use of the Pythagorean theorem. 

The reader will note that the proofs we give for both Theorem A and B rely on a similar argument (coming from \cite{Na1}), namely, putting the 
polygonal piece (the ziggurat or the pyramid) on the hypotenuse toward the interior of the original triangle, and then looking for (scaled) 
rotated versions of the later that naturally arise in the resulting configuration. 

\begin{center}
\includegraphics[width=8cm]{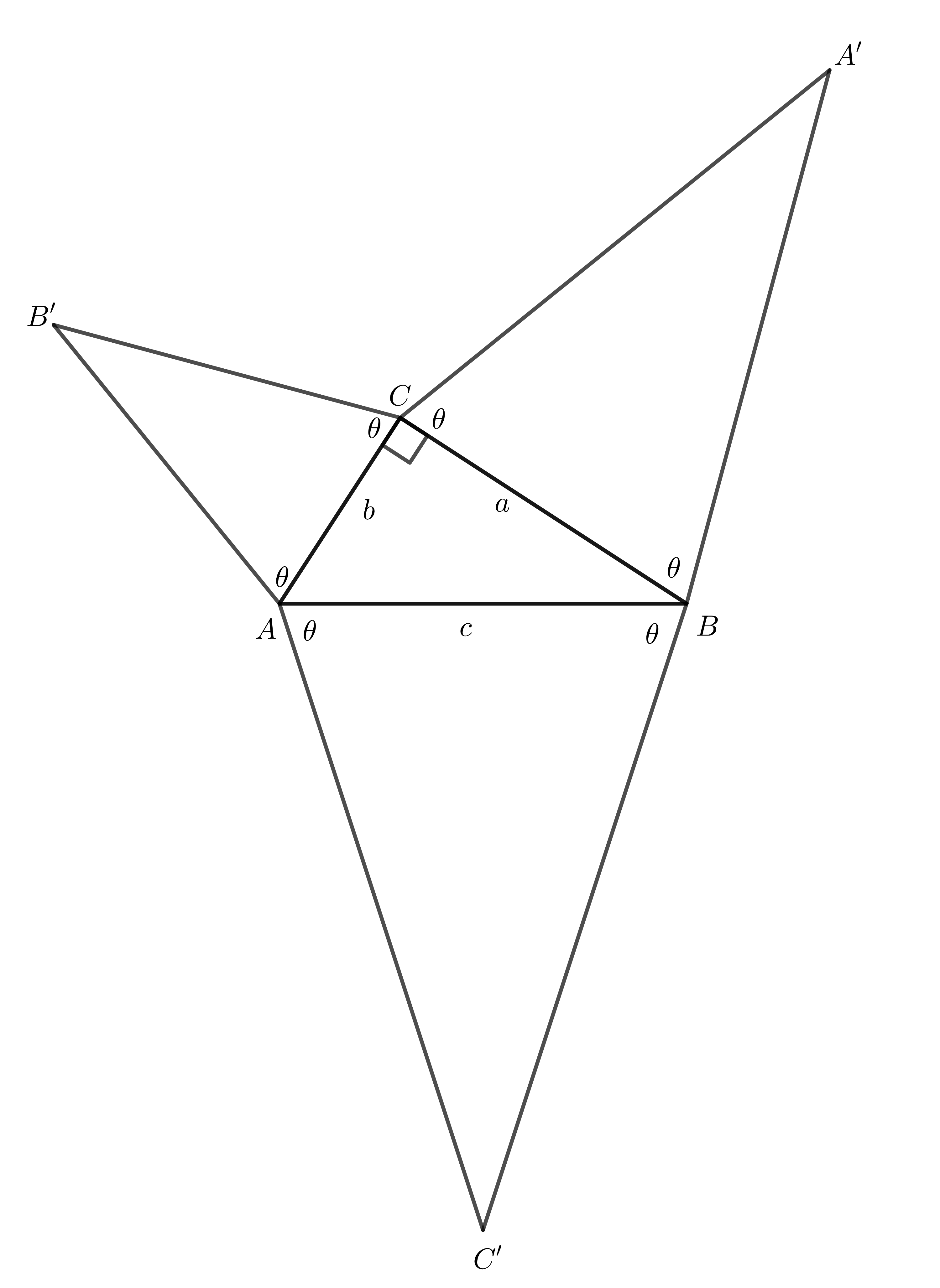}
\end{center}

\section{A general geometric configuration}

Arguably, the main contribution of this note is not to propose a new proof of the Pythagorean theorem (there are so many in the literature!) 
but introducing a new geometric configuration to do this. (See  
\cite{Ma} and references therein on this.) To accomplish such a task here using ziggurats, we 
fix an angle $60^{\mathrm{o}} \leq \theta \leq 135^{\mathrm{o}}$, and we build external ziggurats over the sides $a$ and $b$ and an internal 
ziggurat over the side $c$, all of them with angle $\theta$. Denote these ziggurats as $D' CA D''$, $E' CB E''$ and $F AB G$, as depicted 
below. Also, let $C'$ be such that $E' C'$ has length $b$ and is parallel to $D' C$. Note that $E'C'D'C$ is a parallelogram (which degenerates 
to a segment for $\theta = 135^{\mathrm{o}}$). We will denote $\alpha$ and $\beta$ the angles of $\triangle ABC$ at $A$ and $B$, respectively.
 
\vspace{-0.2cm}
 
\begin{center}
\includegraphics[width=12cm]{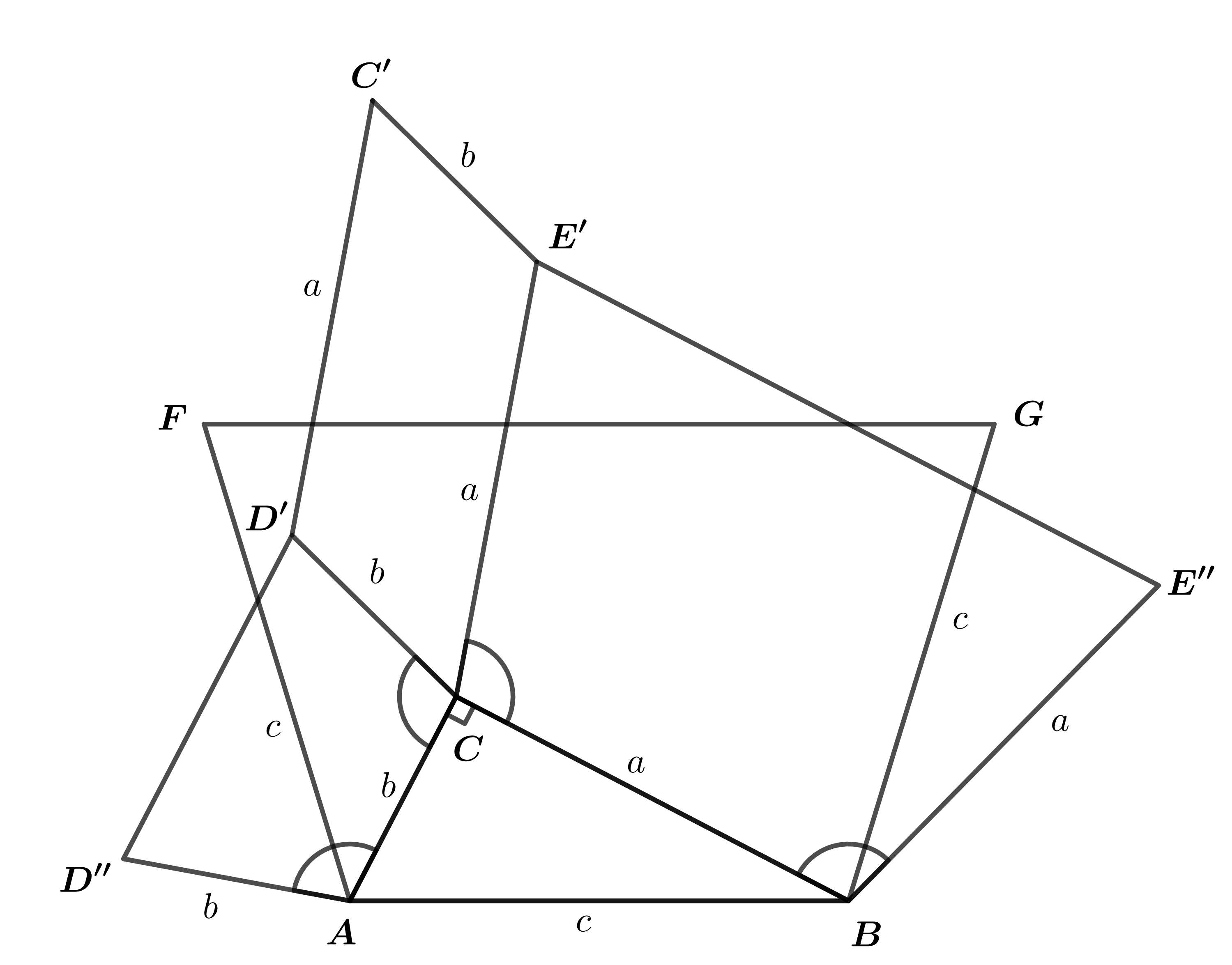}
\end{center} 

 
 \noindent{\bf Lemma.} For $\theta \neq 90^{\mathrm{o}}$, the quadrilaterals $C'D'D''F$ and $C'E'E''G$ are parallelograms 
 (which degenerate to line segments for $\theta = 90^{\mathrm{o}}$).
 
 \vspace{0.35cm}
 
 \noindent{\em Proof.} The triangle $\triangle D''FA$ is a rotated copy of $\triangle ABC$ (it arises from a rotation by an angle $\theta$ about $A$). 
 Therefore, $FD''$ has length $a$,  which is also the length of $C'D'$. We claim that these two segments are also parallel (this implies that 
$C'D'D''F$ is a parallelogram, and the proof for $C'E'E''G$ is analogous). Indeed, this is easily checked by looking at the slopes of these  
segments with respect to $AB$. Namely, on the one hand, $C'D'$ is parallel to the segment $E'C$, whose line is obtained from that of $AB$ 
first by negatively rotating (at $B$) of an angle $\beta$ and later positively rotating (at $C$) of an angle $\theta$. Thus, the slope in this case 
is $\theta - \beta$. On the other hand, the line of $FD''$ is obtained from that of $AB$ first by positively rotating (at $A$) of an angle 
$\theta$ and later negatively rotating (at $F$) of an angle $\beta$. Thus, the slope in this case is again $\theta - \beta$.
$\hfill\square$
 
\vspace{0.18cm} 
 
\begin{center}
\includegraphics[width=10cm]{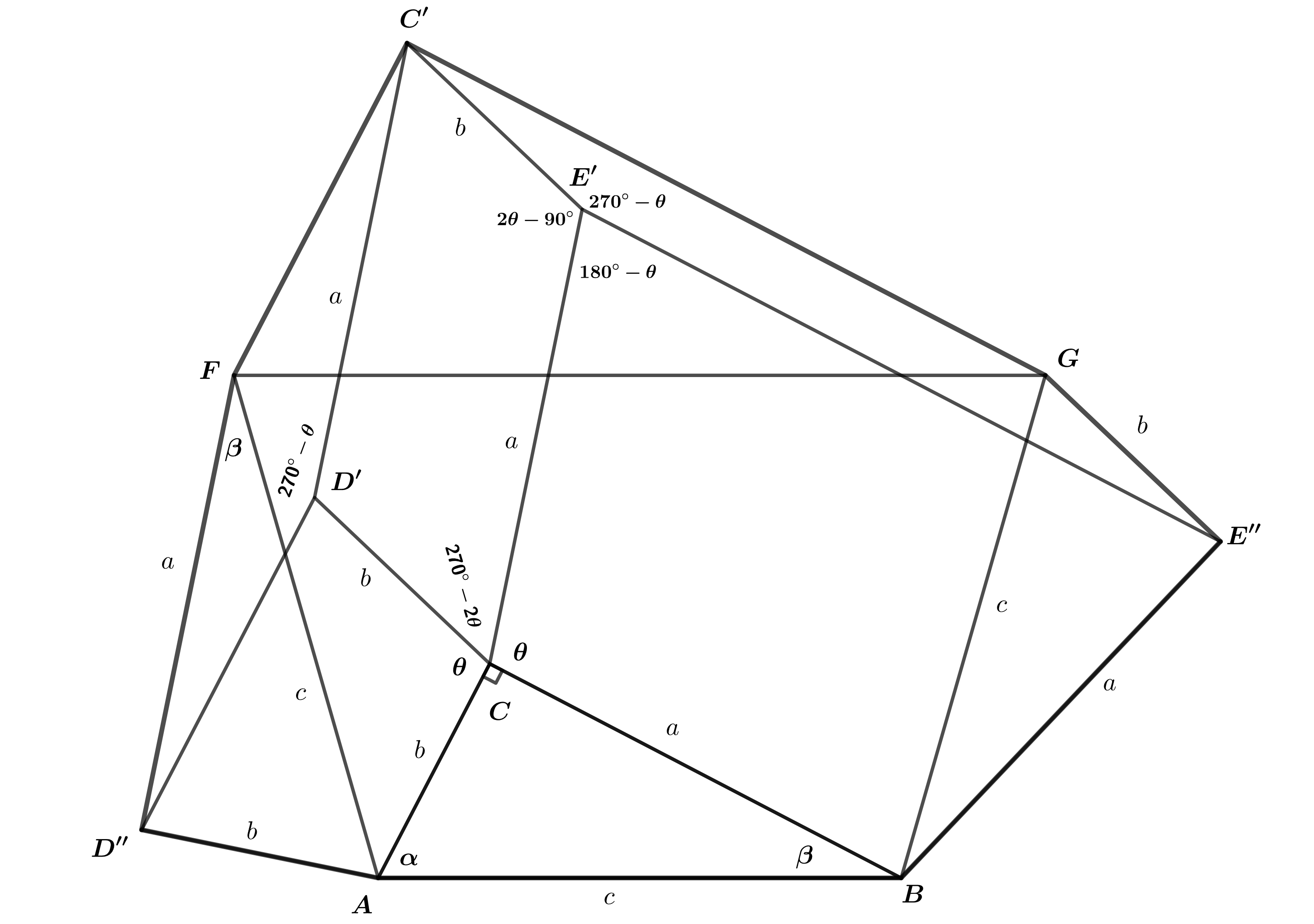}
\end{center}

\vspace{-0.2cm}

\noindent{\bf Observation.} In the proof above, we did not use the hypothesis that the angle at $C$ is a right angle. The 
lemma thus provides a very general configuration. Also, we did not use the restriction on the angle $\theta$, though  
for angles outside the range, the ziggurats may (auto-)intersect. 

\section{The proof argument}

The proof of Theorem A will follow from computing the area of the polygon $ABE''GC'FD''$ (denoted by 
$\mathcal{P}$ in what follows) in two different ways. We will do this in trigonometric terms and later we 
will discuss on the true necessity of using trigonometry.

\vspace{0.2cm}

\noindent{\bf A proof argument for Theorem A.} 
On the one hand, $\mathcal{P}$ consists of the $(c,\theta)$-ziggurat plus the triangles $\triangle GBE'' $, $\triangle AFD''$ 
and $\triangle FGC'$. The first two of these are congruent to $\triangle ABC$, and the latter is similar to it, 
by the lemma above. Moreover, the similarity ratio is easily seen to be equal to $(1-2\cos(\theta))$. Therefore,  
\begin{equation}\label{eq-uno}
area(\mathcal{P}) = area \big( (c,\theta)\mbox{-ziggurat} \big) + \big( 2 + (1-2\cos(\theta))^2 \big) \times area (\triangle ABC).
\end{equation}

On the other hand, $\mathcal{P}$ consists of the $(a,\theta)$ and $(b,\theta)$-ziggurats plus the triangle $\triangle ABC$ and the 
parallelograms $C'D'D''F$, $C'E'E''G$ and $E'C'D'C$. Note that $C'D'D''F$ (resp. $C'E'E''G$) has sides of length $a$ and 
$b (1-2\cos(\theta))$ (resp. $a(1-2\cos (\theta))$ and $b$). Using that $\angle ACB = 90^{\mathrm{o}}$ and chasing angles, 
one easily concludes that 
$$area (E'C'D'C) = ab \, \sin (270^{\mathrm{o}} - \theta)$$
and 
$$area (C'D'D''F) = area (C'E'E''G) = ab \, (1-2 \cos (\theta)) \sin (270^{\mathrm{o}} - \theta).$$ 
Since $ab = 2 \, area (\triangle ABC)$, one obtains that 
\begin{eqnarray*}\label{eq:dos}
area(\mathcal{P}) \!\!\!
&=& \!\!\!area \big( (a,\theta)\mbox{-ziggurat} \big) + 
area \big( (b,\theta)\mbox{-ziggurat}\big) +  \hspace{5.5cm}{(2)} \\
&& \!\!\! + \,\, \big( 1+ 4 \, (1-2\cos(\theta)) \sin (270^{\mathrm{o}} - \theta) 
 + 2 \sin (270^{\mathrm{o}} - 2 \theta)  \big) \times area (\triangle ABC).
 \end{eqnarray*}

A careful but elementary trigonometric manipulation gives\footnote{Note that this is a pure algebraic 
manipulation; in particular, no signed area of polygons is used here.} 
\begin{small}
\begin{eqnarray*}
1+ 4 \, (1-2\cos(\theta)) \sin (270^{\mathrm{o}} - \theta) + 2 \sin (270^{\mathrm{o}} - 2 \theta)
&=& 
1 -  4 (1-2\cos(\theta)) \cos (\theta) - 2 \cos (2 \theta) \\
&=&
1 - 4 \cos (\theta) + 8 \cos^2(\theta)- 2 \, (2 \cos^2(\theta)-1) \\
&=&
3 - 4 \cos (\theta) + 4 \cos^2 (\theta) \\
&=& 
2 + (1-2\cos(\theta))^2. 
\end{eqnarray*}
\end{small}Using this equality to compare (1) and (2), this allows us to conclude that 
$area \big( (c,\theta)\mbox{-ziggurat} \big)$ equals the sum of  
$area \big( (a,\theta)\mbox{-ziggurat} \big)$ plus $area \big( (b,\theta)\mbox{-ziggurat} \big) ,$ 
as claimed.
$\hfill\square$

\vspace{0.34cm}

A short discussion on the use of trigonometry becomes important here, since the reader might be afraid in that 
the manipulations above implicitly use the (trigonometric version of the) Pythagorean theorem. A careful check 
shows that  we used the definition of the trigonometric functions for acute angles, their extensions to arbitrary 
angles, as well as the formulae for the trigonometric functions of sums and differences of angles. 
However, in the last sequence of equalities, we also used the identity 
\begin{equation}\label{eq-problema}
\cos (2 \theta) = 2 \cos^2 (\theta) - 1.
\end{equation}
The problem is that this is based on the Pythagorean identity \, $\sin^2 (\theta) + \cos^2 (\theta) = 1$, 
because the natural formula that arises from the definitions is 
$$\cos (2 \theta) = cos^2 (\theta) - \sin^2 (\theta);$$ 
compare the Remark below.

In order to circumvent this circularity issue, we propose a ``proof (almost) without words'' 
of the identity (\ref{eq-problema}) based on the relation \, $\sin (2 \theta) = 2 \sin (\theta) \cos (\theta)$ \, 
(which is properly available without the use of the Pythagorean theorem; see \cite{Zi}). 
Namely, in the image below, $\triangle ABC$ is similar to $\triangle OB'C'$, and the equality 
$\frac{AB}{CB} = \frac{OB'}{B'C'}$
becomes 
$$1 + \cos (2 \theta) 
= AB 
= CB \cdot \frac{OB'}{B'C'} 
= \sin (2 \theta) \cdot \frac{\cos (\theta)}{\sin (\theta)},$$
from where the identity (\ref{eq-problema}) follows directly.

\vspace{-0.2cm}
\begin{center}
\includegraphics[width=8.4cm]{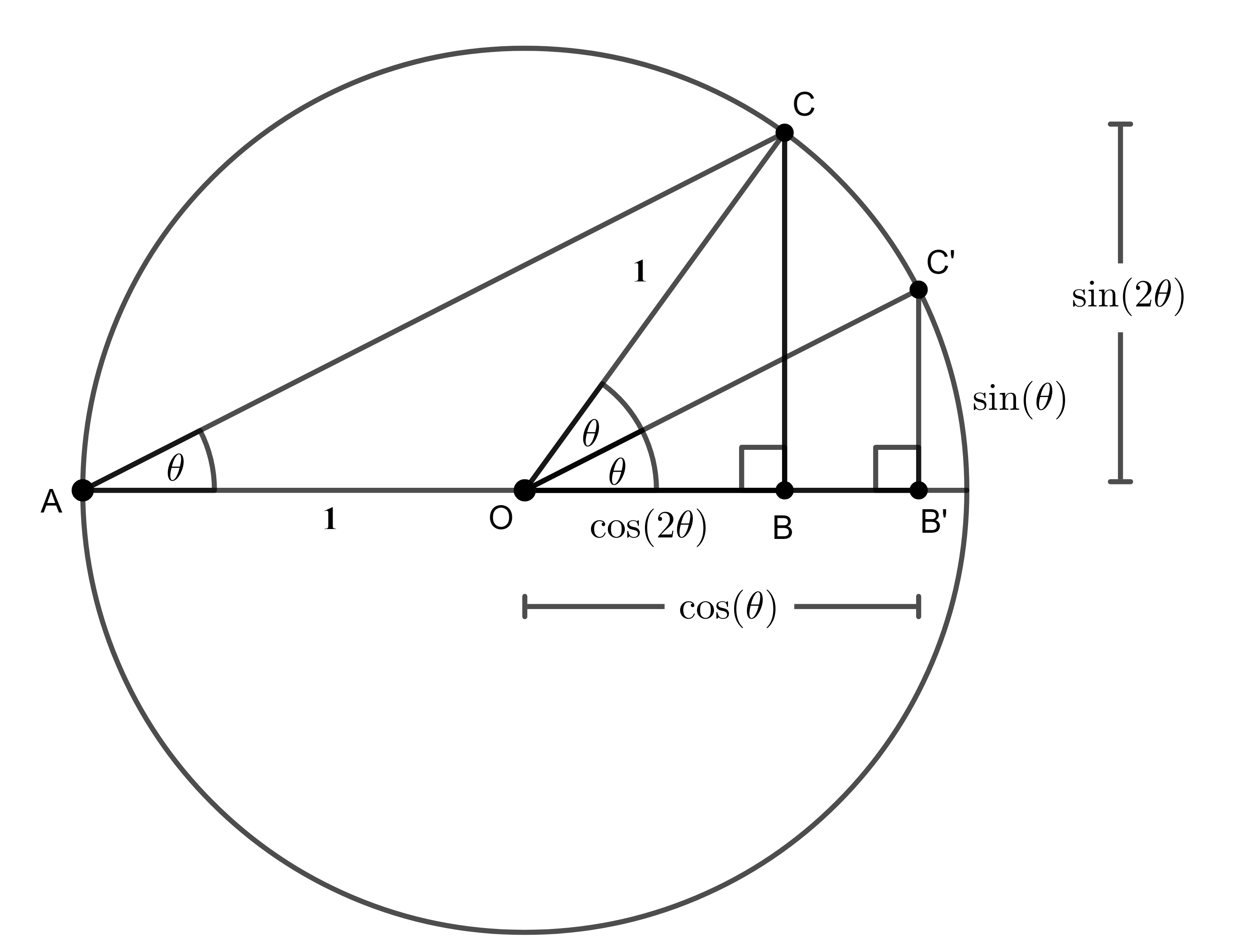}
\end{center}
\vspace{-0.2cm}

The previous argument gives still another trigonometric proof of the Pythagorean theorem, at least for angles no 
larger than $45^{\mathrm{o}}$ (it is not hard to modify the argument for angles between $45^{\mathrm{o}}$ and $90^{\mathrm{o}}$; 
alternatively, one can use the formulae for trigonometric functions of double angles). This is along the lines of \cite{Zi} and  
particularly \cite{Lu}. Since there has been a lot of activity in this direction in recent years (see for instance \cite{JJ}), we 
take this opportunity to include a remark that is somehow unrelated to our geometric construction. 

\vspace{0.25cm}

\noindent {\bf Remark.} If we want to show the Pythagorean theorem for acute triangles using trigonometry, 
a very simple argument would consist in just letting $\alpha' = \alpha$ in the identity 
\begin{equation}\label{eq-stupid}
\cos (\alpha - \alpha') = \cos (\alpha) \cos (\alpha') + \sin (\alpha) \sin (\alpha').
\end{equation}
However, as claimed in \cite{Zi}, this is not allowed, as the trigonometric ratios are not defined for the zero angle... 

Nevertheless, 
nothing prevents us from using a continuity argument, namely by varying $\alpha'$ to make it converge to $\alpha$ from below. 
It is very instructive to look at what happens geometrically when doing this. Indeed, a standard diagram that shows equality 
(\ref{eq-stupid}) is depicted below (see the appendix of \cite{She} for more on this). 

\begin{center}
\includegraphics[width=8.1cm]{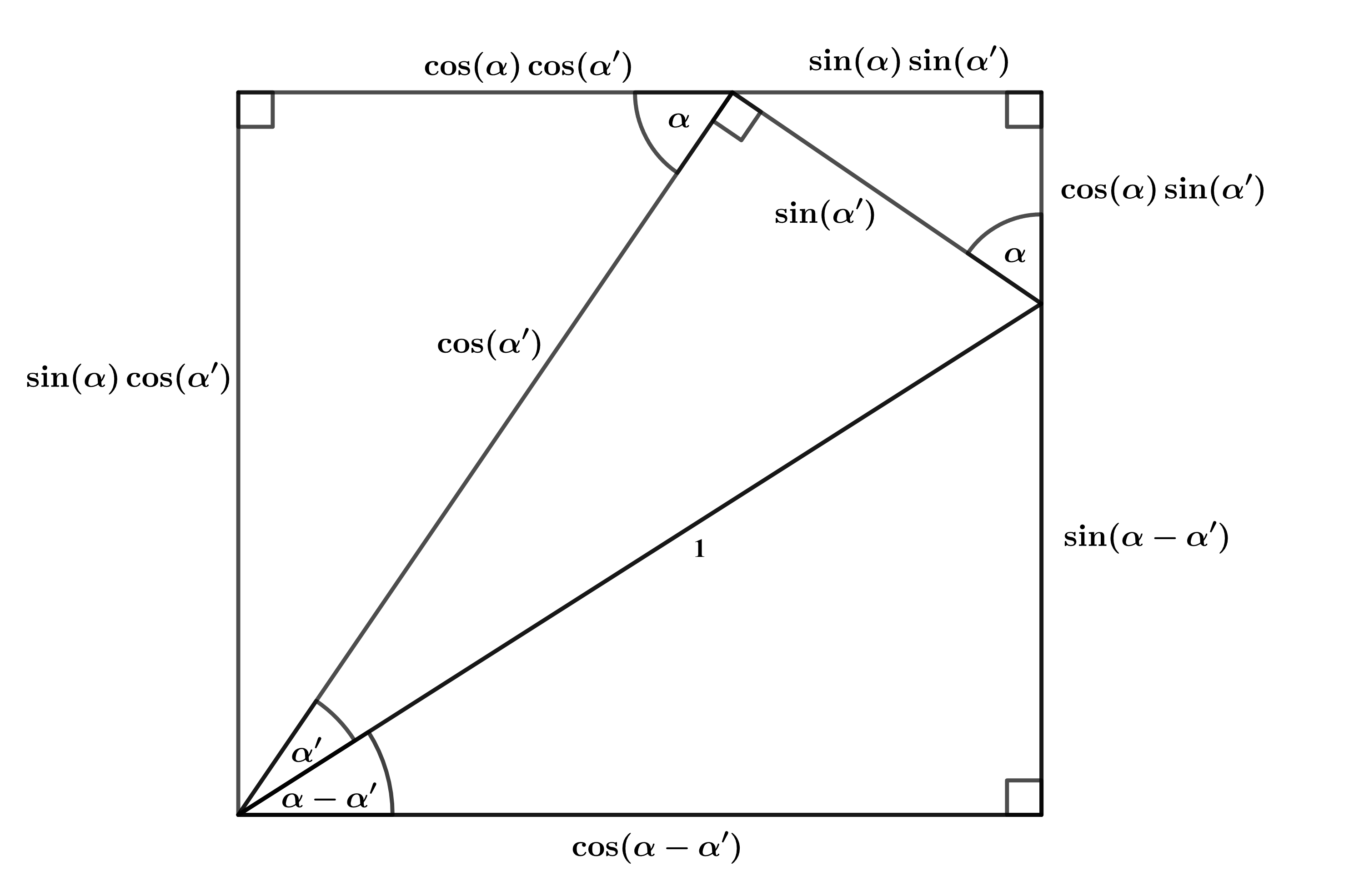} 
\includegraphics[width=8.1cm]{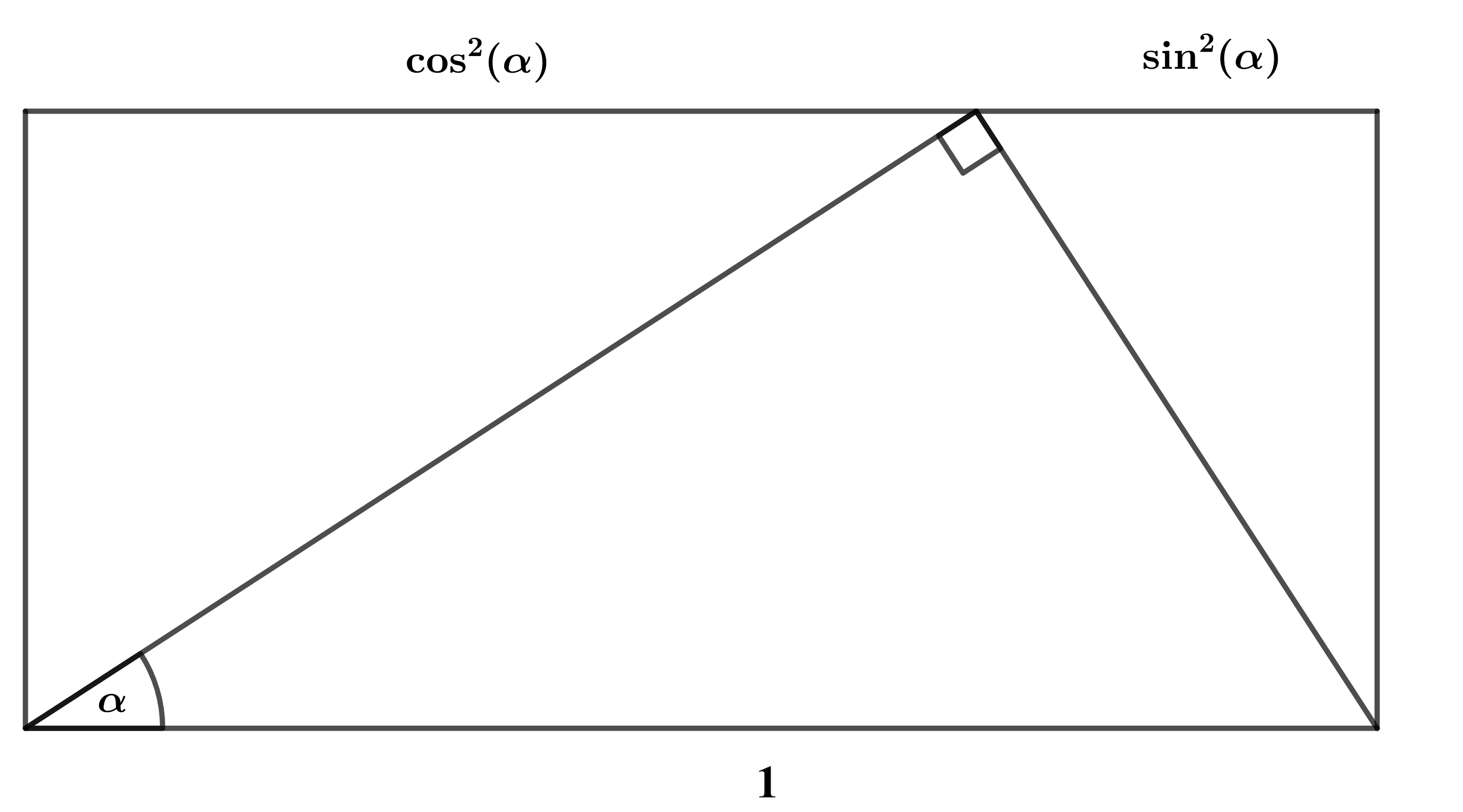}
\end{center}

When passing to the limit, 
this becomes the diagram on the right, from where the Pythagorean relation emerges very clearly once again 
(compare to the Proof 41 in \cite{Bo} as well as \cite{Shi}; see also \cite{Ca} for a more elaborate proof of the 
Pythagorean theorem using analysis).

\section{Pyramids instead of ziggurats}

As for the case of ziggurats, we draw external pyramids over the sides $a$ and $b$ and an internal one over the side $c$.

\vspace{0.25cm}

\noindent{\bf A proof argument for Theorem B.}  
Let us denote by $A'$, $B'$ and $C'$ the vertices of the pyramids, as depicted below. Note that the length of the equal sides of the pyramid 
over $a$ (resp. $b$, $c$) equals $a$ (resp. $b$, $c$) times the constant factor $r = r_{\theta}= \frac{1}{2 \cos (\theta)}$. If we positively rotate 
$\triangle ABC$ about $A$ and then apply a homothety of ratio $r$, we obtain $\triangle AC'B'$. In particular, the length of $B'C'$ equals $ra$.
Analogously, if we negatively rotate $\triangle ABC$ about $B$ and later we apply a homothety of ratio $r$, we obtain $\triangle C'BA'$. Hence, the 
length of $C'A'$ equals $rb$. As a consequence, $CA'C'B'$ has opposite sides of equal length, and therefore it is a parallelogram (which 
degenerates to a segment for a specific value of $\theta$). Note that so far we have not used the fact that the angle at $C$ is a right angle.

Let us now definitively introduce the value $90^{\mathrm{o}}$ for the angle at $C$, and let us chase angles. As before, we can think of 
the polygon $\mathcal{P} = ABA'C'B'$ in two different ways:

\noindent - As the union of $\triangle ABC$, the $(a,\theta)$ and $(b,\theta)$-pyramids, and the parallelogram $CA'C'B'$. As such, 
its area equals the sum of the areas of the two small pyramids plus 
$$\frac{ab}{2} + ra \cdot rb \cdot \sin (2 \theta - 90^{\mathrm{o}}) = ab \left( \frac{1}{2} - r^2 \cos (2 \theta) \right).$$

\noindent - As the union of the $(c,\theta)$-pyramid and the triangles $\triangle AC'B'$ and $\triangle B A' C'$. Viewed this way, its area 
equals that of the $(c,\theta)$-pyramid plus 
$$2 \, \frac{ra \cdot rb}{2} = r^2 ab.$$

\noindent Using that $r = 1 / (2 \cos (\theta)))$, one easily checks that the last two expressions above coincide, thus establishing the result. Note 
that, again, one is forced to use the trigonometric identity $\cos (2\theta) = 2\cos^2 (\theta) - 1$ in this computation. $\hfill\square$

\begin{center}
\includegraphics[width=9cm]{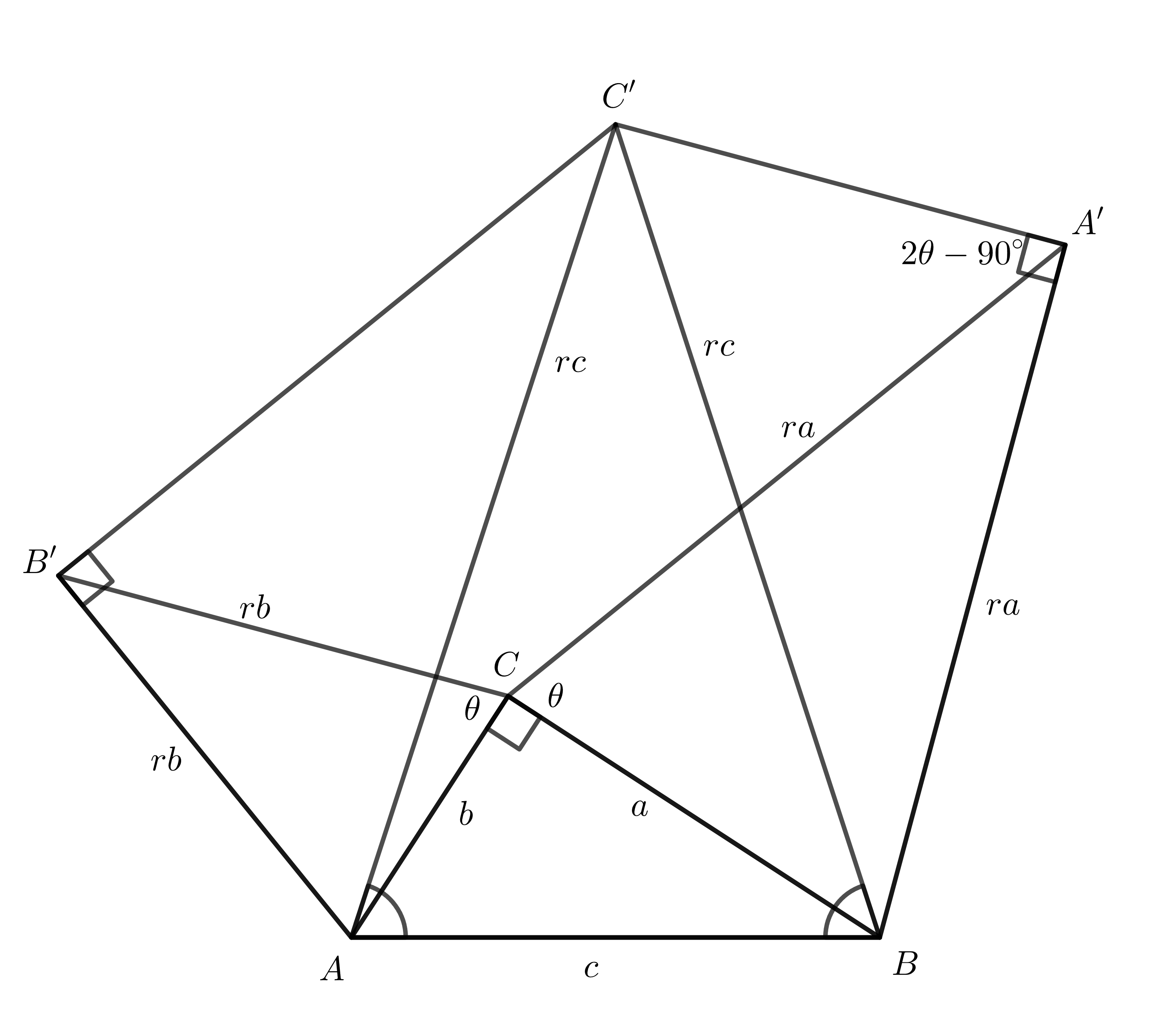}
\end{center} 



\section{Further remarks on the geometric configurations}

To conclude, we go quickly into a deeper analysis of the pictures that arise for specific angles $\theta$. As we will see, we will retrieve 
several very-well known configurations leading to the Phytagorean theorem, some of which fall into the category of ``proofs without words'' 
(we provide the explicit references we found for each case, but certainly there are many others !). The reader can check that in all the cases 
described below, the trigonometric computations involved in the proof become elementary (in particular, the configurations can be drawn 
with a rule and a compass).

\vspace{0.35cm}

\noindent{\bf Configurations via ziggurats.} We start with five configurations arising from the use of ziggurats. 

\vspace{0.25cm}

\noindent{\bf {\underline{$\theta = 90^{\mathrm{o}}$:}}} For this choice, the classical configuration below emerges. 
The proof argument then becomes very close to that of the proofs listed as 24, 63, 69 and 70 in \cite{Bo}.

\begin{center}
\includegraphics[width=7.2cm]{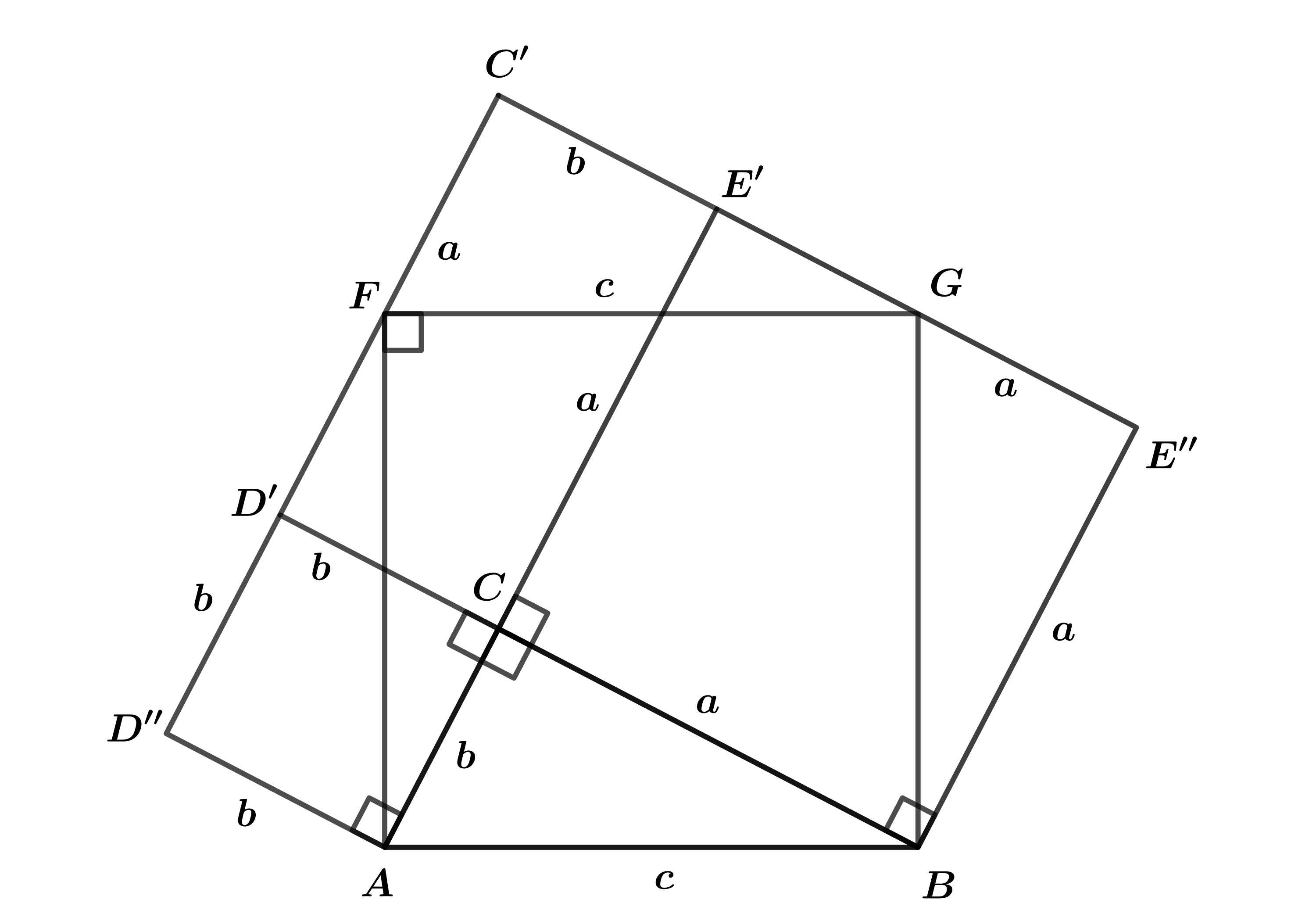}
\end{center}

\noindent{\bf {\underline{$\theta = 60^{\mathrm{o}}$:}}} Here, the configuration that arises and the proof argument 
are those developed in \cite{Na1} (see  \cite{Na2} for an animation). See the picture on the left below. 

\begin{center}
\includegraphics[width=8.2cm]{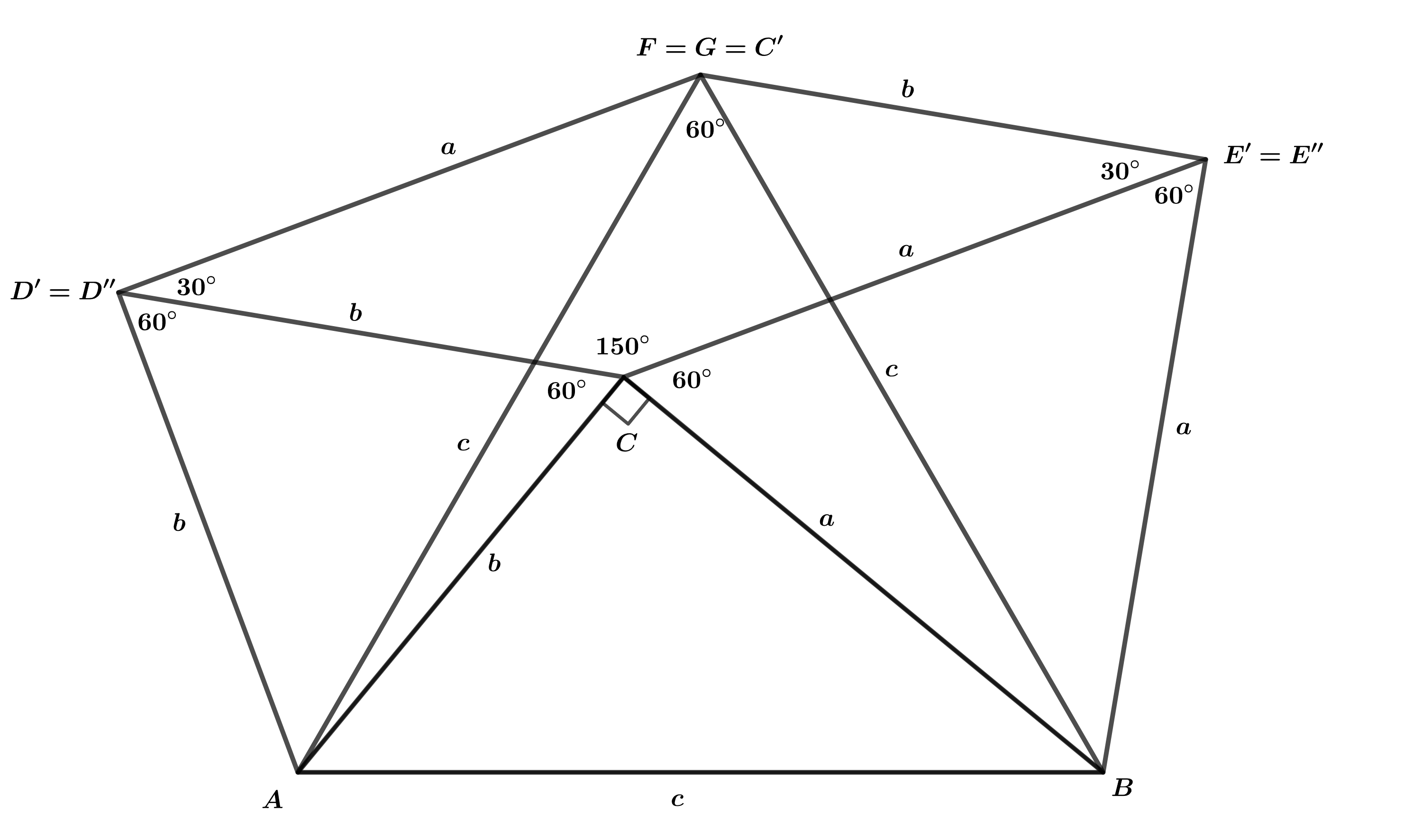}
\includegraphics[width=8cm]{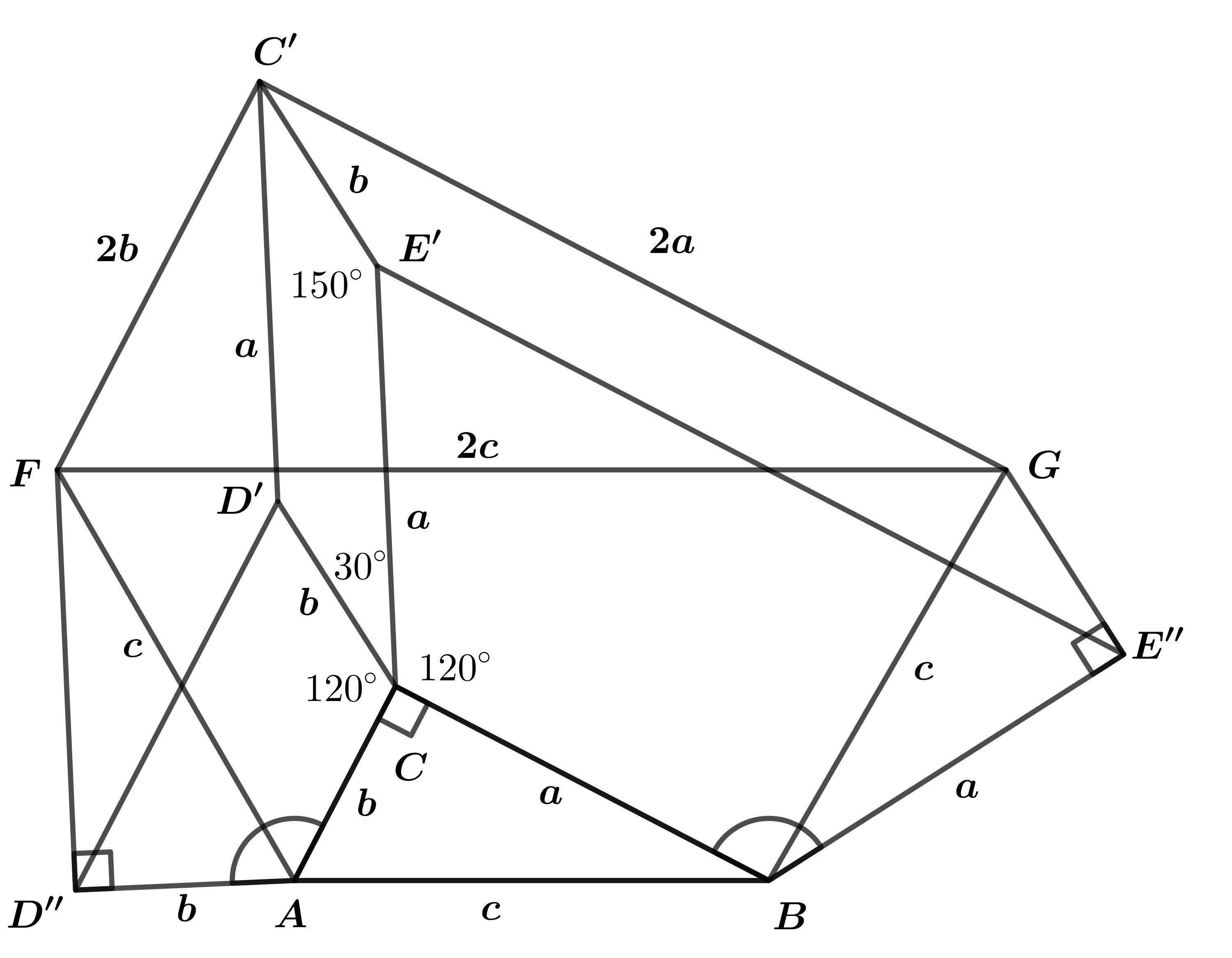}
\end{center}

\noindent{\bf {\underline{$\theta = 120^{\mathrm{o}}$:}}} Each ziggurat corresponds to half of the 
corresponding regular hexagon. Note that the area of the parallelogram $CE'C'D'$ coincides with that of $\triangle ABC$, 
while the areas of each $C'D'D''F$ and $C'E'E''G$  are twice this. See the picture on the right above.

\vspace{0.3cm}

\noindent{\bf {\underline{$\theta = 135^{\mathrm{o}}$:}}} Each ziggurat corresponds to a precise fraction (the ``basis'')  
of the corresponding regular octagon; see the picture on the left below. One can readily check that the equality of areas 
that arose in the proof argument of Theorem A boils down to the elementary algebraic identity \, 
$2 + (1+\sqrt{2})^2 \, = \, 1 + 2 \,  \sqrt{2} \,  (1+\sqrt{2}).$

\begin{center}
\includegraphics[width=8cm]{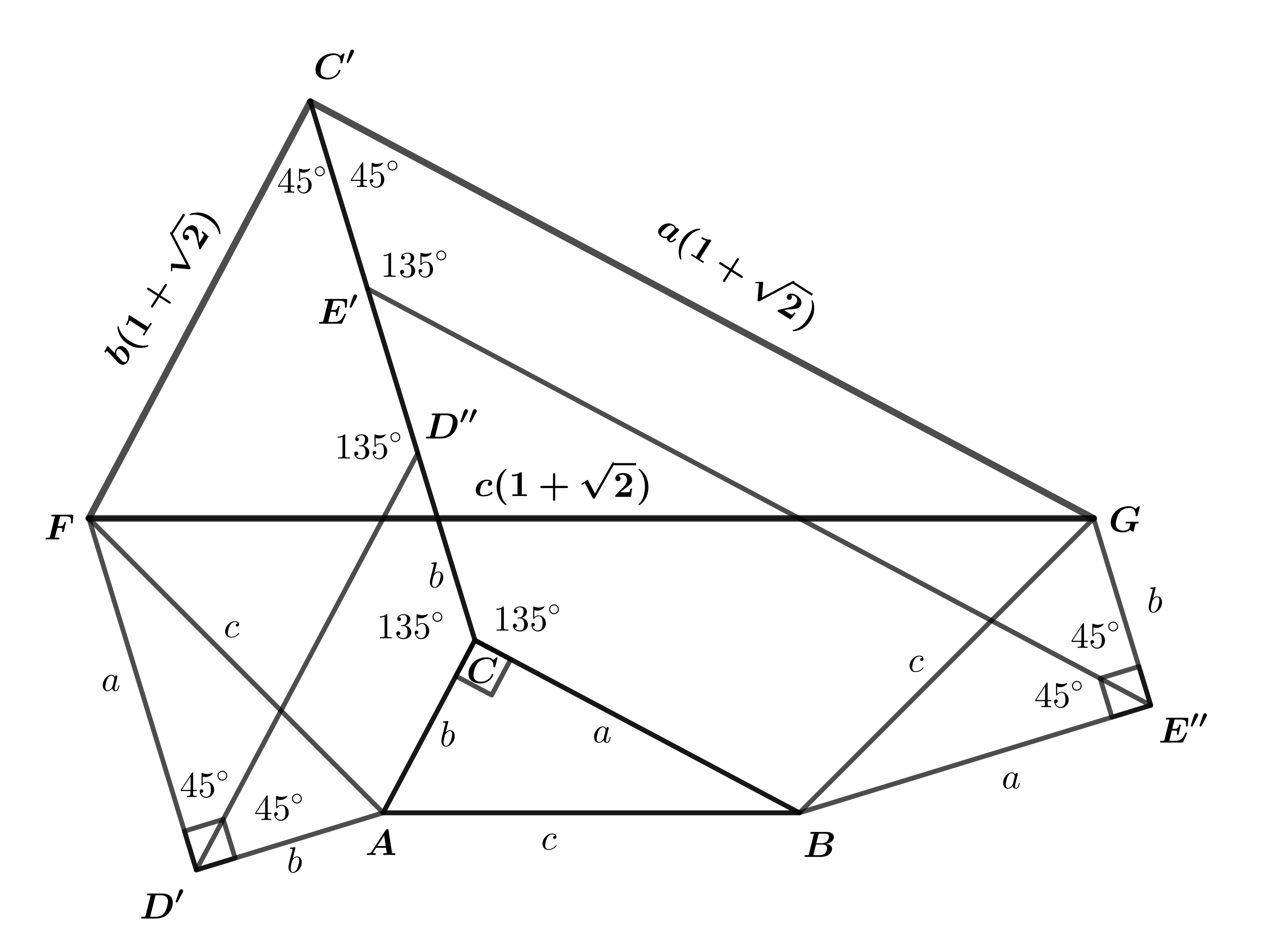}
\includegraphics[width=7.8cm]{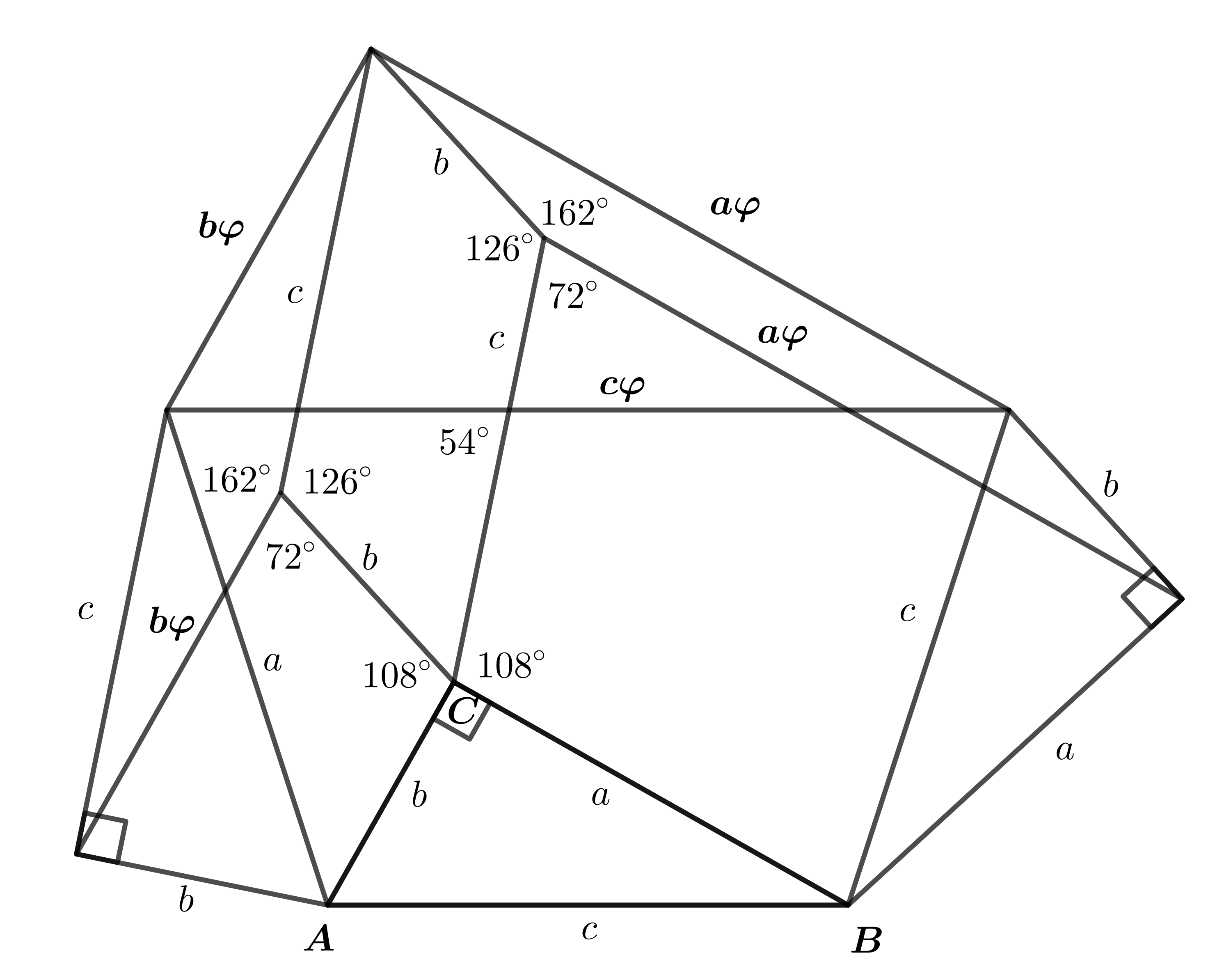}
\end{center}

\vspace{0.3cm}

\noindent{\bf {\underline{$\theta = 108^{\mathrm{o}}$:}}} Each ziggurat corresponds to the ``basis'' of a regular pentagon. 
Having in mind that the length of the diagonal of such a pentagon equals $\varphi$ times the length of its side (where $\varphi$ denotes 
the golden mean $(1+\sqrt{5})/2$), it is straightforward to check (using that $\sin (18^{\mathrm{o}}) = 1/2\varphi$ and 
$\cos (36^{\mathrm{o}}) = \varphi / 2$) that the equality of areas in the proof comes from an 
algebraic relation that ultimately reduces to $\varphi^2 = 1 + \varphi$. 

\vspace{0.32cm}

\noindent{\bf Configurations via pyramids.} As for the case of the ziggurats, there are some configurations arising from the argument 
above that we would like to highlight. 

\vspace{0.2cm}

\noindent{\bf {\underline{$\theta = 45^{\mathrm{o}}$:}}} The configuration that arises goes along the lines 
of the classical ones. Actually, the proof that arises is very close to Proof 64 in \cite{Bo}.

\begin{center}
\includegraphics[width=9cm]{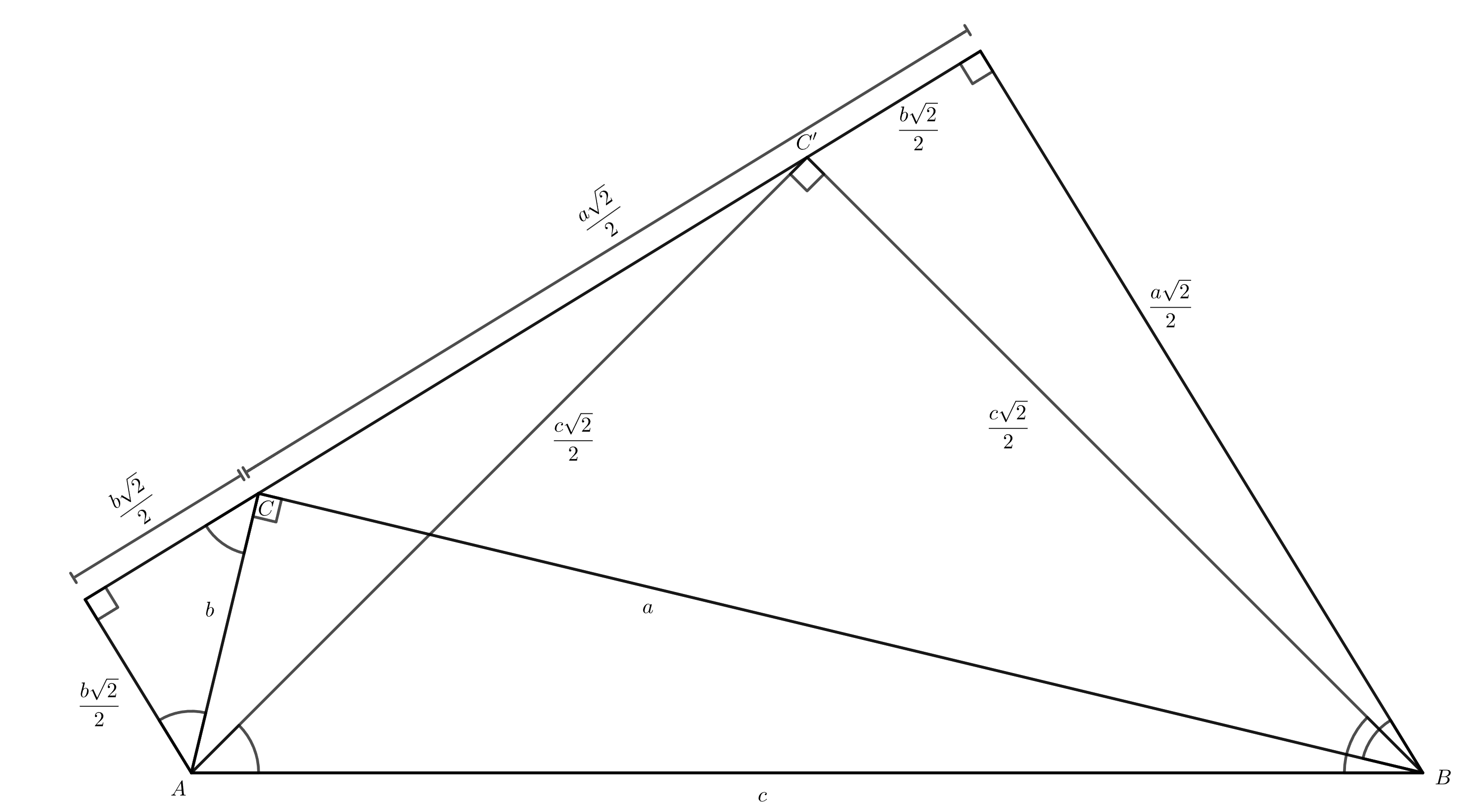}
\end{center} 

\vspace{0.32cm}

\noindent{\bf {\underline{$\theta = 60^{\mathrm{o}}$:}}} Nothing is new here, as an $(\ell,60^{\mathrm{o}})$-pyramid is the 
same as an $(\ell,60^{\mathrm{o}})$-ziggurat. We retrieve again the configuration and the proof argument from \cite{Na1}. 

\vspace{0.32cm}

\noindent{\bf {\underline{$\theta = 72^{\mathrm{o}}$:}}} For this choice, the $(\ell,72^{\mathrm{o}})$-pyramid is a golden triangle (with basis 
$\ell$).  The argument above then ``directly'' establishes that the sum of the areas of the golden triangles over $a$ and $b$ equals the area 
of the golden triangle over $c$. Now, noting that a regular pentagon is the union of two basic ziggurats minus a ``central'' golden triangle, 
this allows establishing the same statement for the regular pentagons over the sides $a,b,c$ in a ``more transparent'' way 
than the one exhibited just using the ziggurats of angle $108^{\mathrm{o}}$. It may be interesting to explore whether this can 
lead to a more natural proof just using a dissection argument (without relying on the Bolyai-Gerwien-Wallace theorem). 
Also, producing a Geogebra applet  to explore this and other configurations arising from this note may be a very pedagogical task.

\vspace{0.25cm}



\vspace{0.3cm}

\noindent{\bf Acknowledgments.} 
I'm indebted to Hugo Caerols, Eduardo Reyes and Michele Triestino for discussions and corrections, 
and to Patrick Popescu-Pampu for his interest and for pointing out to me the reference \cite{Ma}. 
I would like to thank all the school students participating in the ``Entrenamiento Matem\'atico Ol\'impico'' 
at USACH for creating a pleasant atmosphere for thinking about the content of this Note.  
This work was supported by the DICYT Regular Project 042532NF$_{-}$REG ``On the growth of certain 
dynamical cocycles'' (Vicerrector\'ia de Investigaci\'on, Innovaci\'on y Creaci\'on VRIC). In particular, the images 
were produced by Juan Pablo Pincheira, a Master student at USACH. Also, the YouTube video linked to 
in \cite{Na2} was produced (several years ago) by Nicolé Geyssel, at that time also a Master student at 
USACH. I would like to express my warm gratitude to both of them.


\begin{footnotesize}

\vspace{0.1cm}

\noindent Andr\'es Navas \\ 

\noindent Departamento de Matem\'atica y Ciencia de la Computaci\'on\\

\noindent Universidad de Santiago de Chile (USACH)\\ 

\noindent Alameda 3363, Santiago, Chile\\ 

\noindent email: andres.navas@usach.cl

\end{footnotesize}

\end{document}